\documentclass{proc-l}
\usepackage{amsmath,amsthm,amsfonts,amssymb}

\newcommand\RR{{\mathbb R}}

\newcommand\NN{{\mathbb N}}

\def\d={\,:=\,}

\newcommand{\lpraum}[2]
 {{{\rm L}}^{#1}(#2)}

\font\frakten=eufm10
\newfam\frakfam
\textfont\frakfam=\frakten

\newtheorem{thm}{Theorem}[section]
\newtheorem{lemma}[thm]{Lemma}
\newtheorem{cor}[thm]{Corollary}
\newtheorem{prop}[thm]{Proposition}

\theoremstyle{definition}

\theoremstyle{remark}
\newtheorem{ex}[thm]{Example}
\newtheorem{rem}[thm]{Remark}

\begin{document}


\title{Admissible vectors for the regular representation}
\author{Hartmut F\"uhr}
\address{Zentrum Mathematik \\ TU M\"unchen \\ D-80290 M\"unchen \\ Germany}
\email{fuehr@ma.tum.de}

\date{October 5, 2000}

\keywords{Continuous wavelet transforms, coherent states, square-integrable
representations, Plancherel theory, cyclic vectors}

\begin{abstract}
 It is well known that for irreducible, square-integrable representations
 of a locally compact group, there exist so-called admissible vectors
 which allow the construction of generalized continuous wavelet transforms.
 In this paper we discuss when the irreducibility requirement can be
 dropped, using a connection between generalized
 wavelet transforms and Plancherel theory. For unimodular groups with type I
 regular representation, the existence of admissible vectors is
 equivalent to a finite measure condition. The main result
 of this paper states that this restriction disappears in the nonunimodular
 case: Given a nondiscrete, second countable group $G$ with type I regular
 representation $\lambda_G$, we show that $\lambda_G$ itself
 (and hence every subrepresentation
 thereof) has an admissible vector in the sense of wavelet theory
 iff $G$ is nonunimodular.
\end{abstract}

\maketitle

\section*{Introduction}

 This paper deals with the group-theoretic approach to the construction
 of continuous wavelet transforms. Shortly after the continuous wavelet
 transform of univariate functions had been introduced, Grossmann,
 Morlet and Paul established the connection to the representation theory
 of locally compact groups
 \cite{GrMoPa}, which was then used by several authors to construct
 higherdimensional analogues \cite{Mu,Bo,BeTa,Fu,Ji}.

 Let us roughly sketch the general
 group-theoretic formalism for the construction of wavelet transforms.
 Given a unitary representation $\pi$ of the locally compact group $G$
 on the Hilbert space ${\mathcal H}_{\pi}$,
 and given a vector $\eta \in {\mathcal H}_{\pi}$ (the representation
 space of $\pi$), we can study the mapping
 $V_{\eta}$, which maps each $\phi \in {\mathcal H}_{\pi}$ to the
 bounded continuous function $V_{\eta} \phi$ on $G$, defined by
 \[ V_{\eta} \phi (x) := \langle \phi, \pi(x) \eta \rangle ~~.\]
 Whenever this operator $V_{\eta}$ is an isometry of ${\mathcal H}_{\pi}$
 into $\lpraum{2}{G}$, we call it a {\bf continuous wavelet transform},
 and $\eta$ is called wavelet or {\bf admissible vector}. 
 The construction has also been studied in mathematical physics, where
 admissible vectors are known under the name {\bf coherent states}.

 Note that by its very definition a wavelet transform is an intertwining
 operator between $\pi$ and $\lambda_G$, the left regular
 representation.
 The observation made in \cite{GrMoPa} was that admissible vectors
 always exist, when $\pi$ is (unitarily equivalent to) an
 {\em irreducible} subrepresentation of the regular representation of $G$.
 Such representations are usually called ``square-integrable'' or
 ``discrete series representations''.
 This covers both the standard continuous wavelet transform on
 ${\rm L}^2(\RR)$, where the underlying
 group is the affine group of the real line, and the windowed Fourier
 (or Gabor) transform, based on an irreducible representation of the
 Weyl-Heisenberg group.

 While this formalism is flexible enough to allow for a variety of
 transforms (as documented by the above cited higherdimensional analogues), 
 several approaches exist to construct transforms in more general settings:
 On the one hand, the square-integrability requirement can be replaced
 by square-integrability on quotients (see the book \cite{AlAnGa}
 for an exposition of these techniques).

 On the other hand, certain
 nonirreducible representations have been considered as well
 \cite{IsKl,KlSt},
 indicating that the irreducibility requirement can be weakened.
 Another step in this direction is taken in the paper \cite{LWWW},
 which studies in full generality the quasiregular representation of
 a semidirect product $\RR^n \rtimes H$ on ${\rm L}^2(\RR^n)$, where $H$ is
 a closed matrix group, with the aim of establishing admissibility
 conditions. 
 \cite{LWWW} gives an almost complete characterization of the
 matrix groups $H$ for which an admissible vector exists; only
 for a small part of those groups the quasiregular representation
 is in fact irreducible.

 These examples indicate a growing interest in wavelet transforms
 arising from reducible
 representations, and they serve as motivation for this paper.
 We give a complete characerization of the subrepresentations of
 the left regular representation of a locally compact group allowing
 admissible vectors, whenever the regular representation is type I.
 The criteria are given in terms of the Plancherel measure of the group.
 Since we
 are dealing with subrepresentations of the regular representation,
 it seems quite natural to employ the decomposition of the regular
 representation into irreducibles, that is, the Plancherel decomposition.
 The connection was already noted by Carey \cite{Ca}, but 
 had not been further pursued.
 To motivate the approach via Plancherel theory, let us consider
 the following toy example:
 \begin{ex}
 \label{toy_ex}
  Let $G=\RR$, and let ${\mathcal H} \subset {\rm L}^2(\RR)$ be a
  translation-invariant closed subspace. Then there exists a 
  measurable subset $U \subset \RR$ such that
  \[ {\mathcal H} = \{ f \in  {\rm L}^2(\RR): \widehat{f} \mbox{ vanishes
  outside } U \} ~~.\]
  The admissible vectors are easily identified with the aid of the
  Fourier transform: For $\eta, \phi \in {\mathcal H}$ we have
  $V_{\eta} \phi = \phi \ast \tilde{\eta}$, with $\tilde{\eta} (x)
  = \overline{\eta(-x)}$, and hence $\widehat{V_\eta \phi}(\omega)
  = \widehat{\phi}(\omega) \overline{\widehat{\eta}}(\omega)$.
  Thus, $\eta$ is admissible for ${\mathcal H}$ iff $|\widehat{\eta}(\omega)|$
  equals one almost everywhere on $U$, and such vectors exist iff
  $U$ has finite Lebesgue measure.

  Another observation will be useful for the following: Denote
  the restriction of the regular representation to ${\mathcal H}$ by
  $\pi$ and let $\widehat{\pi}$ denote the representation obtained
  by conjugating $\pi$ with the Fourier transform, so that 
  $\widehat{\pi}$ operates on ${\rm L}^2(U)$ by $\widehat{\pi}(x)
  f(\omega) = e^{i\omega x} f(\omega)$. If we choose two vectors
  $\eta, \phi \in {\rm L}^2(U)$, then we see that 
  \[ (V_{\eta} \phi) (x) = \int_U \phi(\omega) \overline{\eta (\omega)}
  e^{-i \omega x} d\omega ~~,\]
  i.e., in this realization of the representation, the wavelet
  transform is just a special instance of Fourier inversion.
 \end{ex}

 It turns out that this argument generalizes almost verbatim to the
 case of unimodular locally compact groups with type I regular
 representation. Just as for $G = \RR$, the Plancherel decomposition
 of such a group allows 
 \begin{itemize}
 \item direct access to all invariant subspaces of ${\rm L}^2(G)$;
 \item to convert convolution operators to pointwise multiplication operators;
 and thus:
 \item to construct admissible vectors on the Plancherel transform side;
 and finally:
 \item to classify the subspaces having admissible vectors by a finite
 measure condition.
 \end{itemize}
 The classification of the subspaces with admissible vectors
 can be found in Theorem \ref{unimod} below.
 The only price we have to pay for passing from abelian to non-abelian
 unimodular groups is that irreducible representations are no longer
 onedimensional. Hence we have to deal with multiplicities in the
 Plancherel decomposition, and to replace scalar multiplication on
 the Fourier side by operator multiplication. 

 However, when we consider nonunimodular groups, the situation changes
 drastically. The complications arise from 
 the family of unbounded operators intervening in the Plancherel transform,
 the so-called {\bf Duflo-Moore} or {\bf formal degree} operators.
 These operators also show up in the decomposition of convolution
 operators. However, it turns out that their unboundedness can be exploited to
 construct admissible vectors for arbitrary subrepresentations of the
 regular representation. This includes the regular representation itself,
 and in fact it is enough to concentrate on this case. These will
 be the main results of this paper:
 \begin{thm}
 \label{main}
  Let $G$ be nonunimodular with type I regular representation 
  $\lambda_G$. Then $\lambda_G$ has an admissible vector. 
 \end{thm}

 \begin{cor}
   Let $G$ be nonunimodular with type I regular representation.
   A representation $\pi$ of $G$ has admissible vectors iff $\pi$ is
   equivalent to a subrepresentation of $\lambda_G$.
 \end{cor}
 \begin{proof}
   Only the ``if''-part remains to be shown.
   So let ${\mathcal H}_{\pi} \subset {\rm L}^2(G)$ be a leftinvariant
   subspace.
   Let $g$ be an admissible vector for the regular representation,
   provided by Theorem \ref{main}.
   If $P$ denotes the projection into ${\mathcal H}_{\pi}$, then $Pg$ is
   an admissible vector for ${\mathcal H}_{\pi}$, since $V_g f = V_{Pg} f$,
   for all $f \in {\mathcal H}_{\pi}$.
 \end{proof}

 Theorem \ref{main} is all the more surprising, as the contrary is proved
 for unimodular groups -- excepting the (trivial)
 case of discrete groups -- with great ease, and without the use of 
 Plancherel theory:
 \begin{prop}
 \label{excl}
  Let $G$ be a unimodular group, such that $\lambda_G$
  has an admissible vector. Then $G$ is discrete.
 \end{prop} 
 \begin{proof}
  Suppose $g$ is an admissible vector. Then $V_g f = f \ast \tilde{g}$,
  and the adjoint of $V_g$ is given by $V_{\tilde{g}}$, where
  $\tilde{g} (x) = \overline{g(x^{-1})}$, which defines an ${\rm L}^2$
  function (here we use unimodularity).
  Hence, for every $f \in {\rm L}^2(G)$, $f = V_{\tilde{g}} (V_{g} f)$,
  and the right-hand side is a continuous function. Hence for every
  ${\rm L}^2$-function $f$ there exists a continuous function which
  coincides with $f$ almost everywhere.
  But then $G$ is discrete.
 \end{proof}

 For unimodular groups with noncompact connected component, we can 
 sharpen the proposition some more: From \cite[Theorem 1.3]{ArLu}, 
 follows that if ${\mathcal H}\subset {\rm L}^2(G)$ has an admissible vector, 
 then ${\mathcal H}$ contains no nontrivial elements supported in a set of
 finite measure.
 
 The proof of the proposition in fact shows that, for any admissible
 vector $g$ for the regular representation of an arbitrary locally compact
 group $G$, $\tilde{g}$ cannot be in ${\rm L}^2(G)$. This indicates that
 the direct construction of admissible vectors will be difficult;
 for one thing, admissible vectors are not compactly supported.

 Let us fix some notation: $G$ denotes a locally compact, second countable
 group. Its left Haar measure is denoted by
 $\mu_G$, $\lpraum{2}{G}$ is the corresponding ${\rm L}^2$-space.
 $\lambda_G$ and $\rho_G$ are the left respectively
 right regular representation; we always assume that $\lambda_G$ is type I.
 $\Delta_G$ denotes the modular function of $G$. Representations
 are always understood to be strongly continuous and unitary.
 $\widehat{G}$ denotes the space of equivalence classes of irreducible
 representations of $G$ endowed with the Mackey Borel structure. We 
 will not explicitly distinguish between representations and equivalence
 classes.

 For a separable Hilbert space ${\mathcal H}$, we let ${\rm dim}({\mathcal H})$
 denote its Hilbert space dimension. ${\mathcal B}_2 ({\mathcal H})$ denotes
 the space of Hilbert-Schmidt operators on ${\mathcal H}$. The norm on
 ${\mathcal B}_2 ({\mathcal H})$ is denoted by $\| \cdot \|_2$. The usual
 operator norm is denoted by $\| \cdot \|_{\infty}$. If an operator
 $A$ is densely defined on ${\mathcal H}$ and has a bounded extension
 to all of ${\mathcal H}$, we denote this extension by $[A]$, and
 we say that ``$[A]$ exists''.

\section{Plancherel Theory}

 In this section we give a short account of Plancherel
 theory, which we then use to reduce the problem of finding admissible
 vectors to the construction of certain operator fields.
 The starting point for the definition of the Plancherel transform is 
 the operator valued Fourier transform on ${\rm L}^1(G)$. Given
 $f \in {\rm L}^1(G)$ and $\sigma \in \widehat{G}$, we define
 \[ {\mathcal F}(f) (\sigma) := \sigma(f) := \int_G f(x) \sigma(x) d\mu_G(x) ~~,\]
 where the integral is taken in the weak operator sense. 
 As direct consequences of the definition we have $\| \sigma(f) \|_{\infty}
 \le \| f \|_1$ and $\sigma(f\ast g) = \sigma(f) \circ \sigma(g)$.

 The Plancherel transform is obtained by extending the Fourier transform
 from ${\rm L}^1(G) \cap {\rm L}^2(G)$ to ${\rm L}^2(G)$. 
 The non-unimodular part of the following Plancherel theorem is
 due to Duflo and Moore \cite[Theorem 5]{DuMo}, whereas the
 unimodular version may be found in \cite{Di}.
\begin{thm}
\label{Pl-Thm}
 Let $G$ be a second countable locally compact group having a type-I
regular representation. Then there exists a measure $\nu_G$ on $\hat{G}$
and a measurable field $(C_{\sigma})_{\sigma \in \hat{G}}$ of selfadjoint
positive operators with densely defined inverse, with the following properties:
\begin{enumerate}
\item[(i)] For $f \in {\rm L}^1(G) \cap {\rm L}^2(G)$ and
 $\nu_G$-almost all $\sigma \in \hat{G}$, the closure of the operator
 $\sigma(f) C_{\sigma}^{-1}$ is a Hilbert-Schmidt operator on
 ${\mathcal H}_{\pi}$. 
\item[(ii)] The map ${\rm L}^1(G) \cap {\rm L}^2(G) \ni f \mapsto
 ([{\sigma(f) C_{\sigma}^{-1}}])_{\sigma \in \widehat{G}}$ extends
 to a unitary equivalence 
 \[ {\mathcal P}: {\rm L}^2(G) \to {\mathcal B}_2^{\oplus} :=
 \int^{\oplus}_{\hat{G}} {\mathcal B}_2({\mathcal H}_{\sigma}) d\nu_G(\sigma) ~~.\]
 This unitary operator is called the {\bf Plancherel transform} of $G$.
 It intertwines the two-sided representation $\lambda_G \times \rho_G$ with
 $\int_{\hat{G}}^{\oplus} \sigma \otimes \overline{\sigma} d\nu_G(\sigma)$.
\item[(iii)] $G$ is unimodular iff almost all $C_{\sigma}$ are scalar
 multiples of the identity operator. In this case we fix $C_{\sigma}
 = {\rm Id}_{{\mathcal H}_{\sigma}}$, and then the measure $\nu_G$ is
 uniquely determined.
\end{enumerate}
\end{thm} 

 In the following $\widehat{f}$ denotes the Plancherel transform 
 of the ${\rm L}^2$-function $f$; in particular in the non-unimodular
 case it should not be confused with the Fourier transform.
 Note that our terminology differs in two ways from \cite{DuMo}: 
 We denote the closure of an operator $A$ by $[A]$, and we use
 slightly different operators. Our $C_{\sigma}$ and the operators
 $K_{\sigma}$ in \cite{DuMo} are related by $C_{\sigma} = K_{\sigma}^{-1/2}$.
 
 In the nonunimodular case, almost every
 $C_{\sigma}$ is an unbounded operator (this will become clear in
 Section 2 below), and it is only fixed up to a constant multiple.
 Since there is apparently no natural choice of normalization, there
 is also no canonical choice of the measure $\nu_G$, which is only
 unique up to equivalence. 

 A further important feature of the Plancherel transform is the decomposition
 of intertwining operators: If $T : {\rm L}^2(G) \to {\rm L}^2(G)$ is a 
 bounded operator which commutes with left translations, then 
 there exists a measurable field of bounded operators $(T_{\sigma})_{\sigma
 \in \widehat{G}}$ with $\| T_{\sigma} \|_{\infty}$ uniformly bounded,
 such that 
 \[ T = \int_{\widehat{G}}^{\oplus} {\rm Id}_{{\mathcal H}_{\sigma}} \otimes
 T_{\sigma} d \nu_G (\sigma) ~~.\]
 This applies in particular to the projection onto invariant subspaces.
 The obvious analogue for the right action of $G$ holds as well.
  
 As explained earlier we are interested in constructing admissible
 vectors on the Plancherel transform side. Hence it is necessary
 to have a pointwise Fourier inversion formula. The following theorem is
 proven in \cite{AlFuKr}. It can be seen as a generalization of 
 \cite[Theorem 4.1]{Li}, or as the nonabelian analogue of the following
 simple fact from
 abelian Fourier analysis: Given an ${\rm L}^2$-function $f$ whose
 Plancherel transform $\hat{f}$ is in ${\rm L}^1$, then, pointwise
 almost everywhere, $f$ equals the Fourier transform of $\widetilde{\hat{f}}$.
\begin{thm}
\label{Pl-Inv}
 Let $A \in {\mathcal B}_2^{\oplus}$ be 
 such that for almost all $\sigma \in \hat{G}$,  $A(\sigma) C_{\sigma}^{-1}$ 
 extends to a trace-class operator. Suppose moreover that the
 mapping $\sigma \mapsto {\rm tr}(|[A(\sigma) C_{\sigma}^{-1}]|)$
 is in ${\rm L}^1(\hat{G},d\nu_G)$.
 Let $a \in {\rm L}^2(G)$ be the inverse Plancherel transform of $A$.
 Then we have (almost everywhere)
\[ a(x) = \int_{\hat{G}} {\rm tr}([A(\sigma)C_{\sigma}^{-1}]
 \sigma(x)^*) d\nu_G(\sigma)~~.
\]
\end{thm}

\begin{rem}
 Let us now briefly explain how wavelet transforms associated to 
 irreducible representations relate to Plancherel transform.
 So let $\pi < \lambda_G$
 be an irreducible subrepresentation, we may assume that $\pi$ is also
 used in the Plancherel deomposition. Then $\nu_G(\{ \pi \}) > 0$, and
 by suitably normalizing the operator $C_{\pi}$ we can assume that
 $\nu_G(\{ \pi \}) = 1$. As is well-known, the admissible vectors
 are those in ${\rm dom}(C_{\pi})$. So let $\eta \in {\rm dom}(C_{\pi})$,
 and associate to $\phi \in {\mathcal H}_{\pi}$ the rank-one operator
 $A_{\pi} = \phi \otimes C_{\pi} \eta$. If we denote by $A \in
 {\mathcal B}_2^{\oplus}$ the operator field which is $A_{\pi}$ at $\pi$
 and zero elsewhere, we find that Theorem \ref{Pl-Inv} is applicable
 and gives for the inverse Plancherel transform $a$ of $A$
 \begin{eqnarray*}  a(x) & = & \int_{\hat{G}} {\rm tr}([A(\sigma) 
 C_{\sigma}^{-1}] \sigma(x)^*) d\nu_G(\sigma) \\
  & = & {\rm tr} ((\phi \otimes C_{\pi} \eta) C_{\pi}^{-1} \pi(x)^*) \\
  & = & {\rm tr} (\phi \otimes \pi(x) \eta) \\
  & = & \langle \phi, \pi(x) \eta \rangle \\
  & = & V_{\eta} \phi (x) ~~.
 \end{eqnarray*}
 Hence $A = (V_{\eta} \phi)^{\wedge}$, and the isometry properties of
 $V_{\eta}$ are easily recognised as special instances of the unitarity
 of Plancherel transform. 

 For multiplicity-free representations, this way of arguing can be
 immediately generalized, and one obtains admissibility conditions which
 are fairly easy to handle. The multiplicity-free case covers for
 instance the quasi-regular representations of semidirect products
 of $\RR^k$ with a closed matrix group,
 acting on ${\rm L}^2(\RR^k)$. In particular the results obtained
 in \cite{IsKl,KlSt} can be interpreted under this perspective.
 
 But also for representations with multiplicities, in particular
 the regular representation itself, one can derive sufficient
 admissibility conditions, which reduce the construction of
 admissible vectors to the construction of certain operator fields.
 This is the subject of the next lemma.
\end{rem}

\begin{lemma}
\label{main_lemma}
 Suppose that $g \in {\rm L}^2(G)$ is such that $[ \hat{g}(\sigma)^*
 C_{\sigma}]$ exists, for almost every $\sigma \in \hat{G}$, and that
 the mapping
 $\sigma \mapsto \| [ \hat{g} (\sigma)^*  C_{\sigma} ] \|_{\infty}$ is
 bounded $\nu_G$-almost everywhere.
 Then, for all $f \in {\rm L}^2(G)$, $V_g f \in {\rm L}^2(G)$, and
 \[ \left( V_g f \right)^{\wedge} (\sigma) = \hat{f}(\sigma) [
  \hat{g}(\sigma)^* C_{\sigma}] ~~.\]
\end{lemma}

\begin{proof}
 For $f \in {\rm L}^2(G)$, consider the measurable field of operators
 $(A(\sigma))_{\sigma \in \hat{G}}$,
 defined by $A(\sigma) =  \hat{f}(\sigma) [\hat{g}(\sigma)^* C_{\sigma}]$.
 Then, by the boundedness condition, $A \in {\mathcal B}_2^{\oplus}$.
 Moreover $[A(\sigma) C_{\sigma}^{-1}] = \hat{f}(\sigma) \hat{g}(\sigma)^*$
 exists and defines a field of trace class operators with 
 integrable trace class norm. Hence we may apply the inversion formula
 and obtain for the inverse Plancherel transform $a$ of $A$:
 \[ a(x) = \int_{\hat{G}} {\rm tr} (\hat{f}(\sigma) \hat{g}(\sigma)^*
 \sigma(x)^*) d\nu_G(\sigma) = \langle f, \lambda(x) g \rangle ~~.\]
\end{proof}

Lemma \ref{main_lemma} gives a sufficient criterion for admissibility:
We have to construct 
$(A(\sigma))_{\sigma \in \hat{G}} \in {\mathcal B}_2^{\oplus}$
in such a way that $B \mapsto B [A(\sigma)^* C_{\sigma}]$
defines an isometry in ${\mathcal B}_2({\mathcal H}_{\sigma})$, for almost
every $\sigma \in \hat{G}$. The latter is easily seen to be equivalent
to the fact that $[A(\sigma)^* C_{\sigma}]^*$ is an isometry.
The inverse Plancherel transform of such an operator field is 
then the desired admissible vector. At this stage, the effect
of nonunimodularity becomes visible: In this case there is a chance
that such Hilbert-Schmidt operators $A(\sigma)$ exist, because
$C_{\sigma}$ is unbounded. In the unimodular case the $C_{\sigma}$
disappear, and thus $A(\sigma)$ has to be both Hilbert-Schmidt (in
particular compact) 
and an isometry, which only works if ${\mathcal H}_{\sigma}$ is finite-dimensional.

Now let us consider the unimodular case. Before we prove the criterion
which generalizes the toy example from above, we need a converse
of Lemma \ref{main_lemma}. The statement is very intuitive and
most likely well-known, but it seems to us that a rigorous proof 
is necessarily somewhat technical, due to the fact that the
Plancherel transform of an arbitrary ${\rm L}^2$-function is
known only almost everywhere, unlike the Fourier transform.
 
\begin{lemma}
\label{main_lemma_converse}
 Let $G$ be unimodular. If $g \in {\rm L}^2(G)$ is such that
 $V_g$ is a bounded operator on ${\rm L}^2(G)$, then $\| \widehat{g} (\sigma)
 \|_{\infty}$ is bounded $\nu_G$-almost everywhere, and we have
 for all $f \in {\rm L}^2(G)$
 \[ (V_g f)^{\wedge}(\sigma) = \widehat{f}(\sigma) \widehat{g}(\sigma)^* 
 ~~ (\nu_G-\mbox{ almost everywhere}) ~~.\]
\end{lemma}

\begin{proof}
 The operator $V_g$ commutes with left translation and hence has
 a decomposition 
 \[ V_g = \int_{\widehat{G}}^{\oplus} {\rm Id}_{{\mathcal H}
 _{\sigma}} \otimes A_{\sigma} d\nu_G(\sigma) ~~,\]
 with a field $(A_{\sigma})_{\sigma \in \widehat{G}}$ of essentially
 uniformly bounded operators. Of course the aim is to show that
 $A_{\sigma} = \widehat{g}(\sigma)^*$ $\nu_G$-almost everywhere.
 For this purpose let $(f_n)_{n \in \NN} \subset {\rm L}^1(G) \cap
 {\rm L}^2(G)$ be a sequence with dense span in ${\rm L}^2(G)$.
 Then, since the Plancherel transform also intertwines the representations
 of the convolution algebra ${\rm L}^1(G)$ arising from the left
 action of $G$, we have for all $n \in \NN$ that  
 \[ \widehat{f_n}(\sigma) A_{\sigma} = (V_g f_n)^{\wedge}(\sigma) =
 (f_n \ast \tilde{g})^{\wedge}(\sigma) = 
  \sigma(f_n) \widehat{g}(\sigma)^* = \widehat{f_n}(\sigma)
 \widehat{g}(\sigma)^*
 ~~,\]
 for all $\sigma$ belonging to a {\it common} conull subset $\Sigma \subset 
 \widehat{G}$. 
 Moreover, possibly after passing to a suitable conull subset of $\Sigma$,
 we can assume that all $A_{\sigma}$ and $\widehat{g}(\sigma)$
 are bounded, and  
 that the span of $\{ \widehat{f_n}(\sigma) : n \in \NN \}$ is
 dense in ${\mathcal B}_2({\mathcal H}_{\sigma})$, for every $\sigma \in
 \Sigma$. 
 But then by the continuity
 of the operators it follows that $A_{\sigma} = \widehat{g}(\sigma)^*$
 for all $\sigma \in \Sigma$. 
\end{proof}

Now we can easily prove the
characterization of the subrepresentations of $\lambda_G$ with
admissible vectors, when $G$ is unimodular. The statement and its
proof are somewhat similar to \cite[Proposition 1.1]{ArLu}; 
compare also \cite[Theorem 2.10]{Ca}.
\begin{thm}
\label{unimod}
 Let $G$ be unimodular.
 Let ${\mathcal H} \subset {\rm L}^2(G)$ be a leftinvariant closed
 subspace, and let $P$ denote the orthogonal projection onto ${\mathcal H}$.
 Then 
 \[  P = \int_{\hat{G}}^{\oplus} {\rm Id}_{{\mathcal H}_{\sigma}} 
  \otimes P_{\sigma} d\nu_G(\sigma) ~~,\]
 with a measurable family of orthogonal projections $(P_{\sigma})_{\sigma 
 \in \hat{G}}$. Then ${\mathcal H}$ has admissible vectors iff 
  almost all $P_{\sigma}$ have finite rank and
 \[ \nu_{\mathcal H} = \int_{\hat{G}} {\rm dim} (P_{\sigma}({\mathcal H}_{\sigma}))
  d\nu_G(\sigma) < \infty ~~. \]
 Every admissible vector $g \in {\mathcal H}$ fulfills $\| g \|^2 =
 \nu_{\mathcal H}$. 
\end{thm}
 
\begin{proof}
 For the sufficiency we note that $(P_{\sigma})_{\sigma \in \hat{G}}
 \in {\rm B}_2^{\oplus}$, and we let $g$ be the inverse Plancherel
 transform of that. Then Lemma \ref{main_lemma} shows that
 $P = V_g$, which means that $V_g$ is the identity operator on
 ${\mathcal H}$, and $g$ is admissible.

 Now let $g$ be an admissible vector for ${\mathcal H}$, and
 define $h = \tilde{g} \ast g$.
 Then we have $P = (V_g)^{\ast} \circ (V_g)
 = V_{\tilde{g}} \circ V_g = V_h$ and $h = V_{\tilde{g}} \tilde{g}
 \in {\rm L}^2(G)$ (note that $V_{\tilde{g}}$ is a bounded operator
 on all of ${\rm L}^2(G)$). Applying Lemma \ref{main_lemma_converse} first 
 to $V_h$ and then to $V_{\tilde{g}}$ yields
 \[ P_{\sigma} = \widehat{h}(\sigma)^* = \widehat{g}(\sigma)^*
 \widehat{g}(\sigma)~~, \]
 $\nu_G$-almost everywhere. Hence
 \begin{eqnarray*}
 \| g \|^2 & = & \int_{\hat G} {\rm tr}(\hat{g}(\sigma) \hat{g}(\sigma)^*)
  d\nu_G(\sigma) \\ & = &
  \int_{\hat G} {\rm tr}(\hat{g}(\sigma)^* \hat{g}(\sigma))
  d\nu_G(\sigma) \\ & = &
  \int_{\hat G} {\rm dim}(P_{\sigma}({\mathcal H}_{\sigma}))  d\nu_G(\sigma) ~~.
 \end{eqnarray*} 
\end{proof}

The number $\nu_{\mathcal H}$ is the ``total Plancherel measure''
associated to the representation $\lambda_G|_{\mathcal H}$. This terminology
expresses the fact that in calculating $\nu_{\mathcal H}$,
each representation occurring in the direct integral decomposition
of $\lambda_G|_{\mathcal H}$ into irreducibles is weighted by its
multiplicity, and then integrated against Plancherel measure.

Combining Proposition \ref{excl} and Theorem \ref{unimod}, we find that 
for nondiscrete unimodular groups the total Plancherel measure of
the regular representation itself is
always infinite. We expect that this has already been noted
elsewhere; it can also be derived from the results 
in the paper by Ludwig and Arnal \cite{ArLu}, where
the quantity $\nu_{\mathcal H}$ is studied in the
context of the so-called qualitative uncertainty property.

\section{Proof of Theorem 0.2}

Throughout this section $G$ is nonunimodular.

We shall now construct the operator field $(\hat{g}(\sigma))_{\sigma \in
\hat{G}}$ with the sufficient properties derived from of Lemma 
\ref{main_lemma}. For this purpose a more detailed knowledge of the
operators $C_{\sigma}$ and the representations involved is indispensable,
and it can be obtained by considering the closed normal, unimodular subgroup
$G_0 := {\rm Ker} (\Delta_G)$. We denote the quotient
$G/G_0$ by $\Gamma$ and endow it with the quotient topology.
It is thus algebraically (but usually not topologically) isomorphic
to some subgroup of $(\RR^+,\cdot)$. By abuse of notation, the
canonical embedding of $\Gamma$ in $\RR^+$ is also denoted by
$\Delta_G$.

The main idea pursued in Paragraph 6 of \cite{DuMo} is to perform a 
Mackey analysis of the group extension $G$ of $G_0$. (A similar
approach was taken by Tatsuuma \cite{Ta}.) 
This means that information on the orbit
space of the natural action of $G$ (actually, $\Gamma$) on the
dual $\widehat{G_0}$ is utilized to derive the Plancherel theory
of $G$ from that of the unimodular subgroup $G_0$.

This approach results in a fairly explicit description of the objects
involved. More precisely,
there exists a $\Gamma$-invariant Borel subset $U \subset \widehat{G_0}$,
which has the following properties:

\begin{enumerate}
\item[(i)] It is $\nu_{G_0}$-conull in $\widehat{G_0}$;
 here $\nu_{G_0}$ is the Plancherel measure of $G_0$.
 The quotient space $U/\Gamma$ is a standard Borel space.
 (This follows combining \cite[Theorem 6, 1.]{DuMo}, \cite[Lemma 13]{DuMo}
 and \cite[Corollary 1 to Theorem 6]{DuMo}.)
\item[(ii)] For every $\sigma \in U$, ${\rm Ind}_{G_0}^G \sigma \in
 \widehat{G}$ \cite[Theorem 6, 1.]{DuMo}.
 Hence $\Gamma$ operates freely on the orbit
 $\Gamma \sigma$ \cite[Lemma 7]{DuMo}.
\item[(iii)] Define $V :=  \{{\rm Ind}_{G_0}^G \rho : \rho \in U \}$,
 then $V$ is a Borel subset of $\hat{G}$ and standard \cite[Theorem 6]{DuMo}.
 By Mackey's theory $V$ is canonically identified with
 the orbit space $U/\Gamma$, and this identification is a Borel
 isomorphism \cite[Lemma 13]{DuMo}. Moreover there exists a Borel
 cross section
 $\tau : V \to U$, such that $\Gamma \times V \ni (\gamma,\sigma)
 \mapsto \gamma \tau(\sigma) \in U$ is a Borel isomorphism
 \cite[Proof of Prop. 10]{DuMo}. In particular $U$ is standard as
 well, and hence $\tau(V)$ is standard \cite[Theorem 3.2]{Ma}.
 
 Using the Borel isomorphism, we transfer the mapping $(\gamma, \sigma)
 \mapsto \Delta_G(\gamma)^{-1/2}$ to $U$ and thus obtain a measurable
 function $\psi$ on $U$. Then $\psi$ obviously fulfills
 $\psi(\gamma \rho) = \Delta_G(\gamma)^{-1/2} \psi(\rho)$, for all
 $\rho \in U$.
\item[(iv)] The Plancherel measure is supported by ${\rm Ind}(U)$,
 and can be obtained from the Plancherel measure of $G_0$ by
 measure disintegration.

 More precisely: Plancherel measure on $V \cong U/\Gamma$ is
 obtained from a measure disintegration along $\Gamma$-orbits,
 as pointed out in \cite[Theorem 6, 3.]{DuMo}: Since $\Gamma$
 operates freely, each orbit $\Gamma \tau(\sigma)$ can be endowed with the 
 image $\mu_{\Gamma \tau(\sigma)}$ of Haar measure under the projection
 map $\gamma \mapsto \gamma \tau(\sigma)$.
 Then there exists a unique measure $\nu_{G}$ on $V$ such that
 \[ d\nu_{G_0} (\rho) = \psi(\rho)^{-2} d\mu_{\Gamma \tau(\sigma)} (\rho)
 d\nu_{G} (\sigma) ~~;\]
 the thus found measure $\nu_{G}$ is the Plancherel measure of $G$.

 Moreover, the operators $C_{\sigma}$ can be given explicitly
 in terms of $\Delta_G$. Since $\Gamma$ operates freely, we
 may realize  $\sigma = {\rm Ind}_{G_0}^G \tau(\sigma)$ on
 $\lpraum{2}{\Gamma, d\mu_{\Gamma}; {\mathcal H}_{\tau(\sigma)}}$, and then
 $C_{\sigma}$ is given by multiplication with $\Delta_G^{-1/2}$:
 \[ (C_{\sigma} \eta) (\gamma) = \Delta_G(\gamma)^{-1/2} \eta(\gamma) ~~,\]
 and ${\rm dom}(C_{\sigma})$ is the set of all $\eta \in \lpraum{2}{\Gamma,
 d\mu_{\Gamma};
 {\mathcal H}_{\tau(\sigma)}}$, for which this product is also in
 $\lpraum{2}{\Gamma, d\mu_{\Gamma}; {\mathcal H}_{\tau(\sigma)}}$. 

 At this point it is easy to see that in the nonunimodular case indeed almost
 every $C_{\sigma}$ is unbounded.
\end{enumerate}

Now let us construct the operator field. We first give
the $A(\sigma)$ pointwise and postpone the questions of measurability
and square-integrability. Given $\sigma \in \hat{G}$, $A(\sigma)$ lives on 
$ {\mathcal H}_{\sigma} = \lpraum{2}{\Gamma, d\mu_{\Gamma};
{\mathcal H}_{\tau(\sigma)}} $.
Pick $c > 1$ in such a way that
$\{ \gamma \in \Gamma: 1 \le \Delta_G^{-1/2}(\gamma) < c \}$ has positive
Haar measure, and define, for $n \in \NN$, $S_n :=
\{ \gamma \in \Gamma: c^n \le \Delta_G^{-1/2}(\gamma) < c^{n+1} \}$.
Let $(u_n^{\sigma})_{n \in \NN} \subset
 \lpraum{2}{\Gamma, d\mu_{\Gamma}; {\mathcal H}_{\tau(\sigma)}} $
be an orthonormal basis. Moreover let
$(v_n^{\sigma})_{n \in \NN} \subset \lpraum{2}{\Gamma, d\mu_{\Gamma}; 
{\mathcal H}_{\tau(\sigma)}}$
be a sequence of unit vectors with ${\rm supp}(v_n^{\sigma}) \subset S_n$.
Define the linear operator $A(\sigma)$ by
\[ A(\sigma)( u_n^{\sigma}) := \| \Delta_G^{-1/2} v_n^{\sigma} \|^{-1} v_n^{\sigma} ~~.\]
Then, since $\| \Delta_G^{-1/2} v_n^{\sigma} \| \ge c^n$, $A(\sigma)$ is a Hilbert-Schmidt
operator. In addition, $C_{\sigma} A(\sigma)$ is an isometry,
since $A(\sigma)$ maps ${\mathcal H}_{\sigma}$ into ${\rm dom}(C_{\sigma})$, and 
$C_{\sigma} A(\sigma)$ maps the orthonormal basis $(u_n^{\sigma})_{n \in \NN}$ to the
orthonormal system $(v_n^{\sigma})_{n \in \NN}$. Finally,
$A(\sigma)^* C_{\sigma}$ has a continuous extension: To see this, let
$f \in {\rm dom} (C_{\sigma})$. $f$ decomposes in the
orthogonal sum $f = g + \sum_{n \in \NN} f_n$,
where $f_n$ has support in $S_n$, and $g$ vanishes on the $S_n$.
By assumption, $\Delta_G^{-1/2} f
\in \lpraum{2}{\Gamma, d\mu_{\Gamma}; {\mathcal H}_{\tau(\sigma)}}$. 
The operator $A(\sigma)^*$ is given by
\begin{eqnarray*}
 A(\sigma)^* (v_n^{\sigma}) & = & \| \Delta_G^{-1/2} v_n^{\sigma} \|^{-1} u_n^{\sigma}\\
 A(\sigma)^* (w)   & = & 0 ~~,~~\mbox{whenever} ~ w \bot \{ v_n^{\sigma} : n \in \NN \}
\end{eqnarray*}
Hence 
\[ A(\sigma)^* C_{\sigma} f = \sum_{n \in \NN} \langle \Delta_G^{-1/2} f_n, v_n^{\sigma} \rangle 
 \| \Delta_G^{-1/2} v_n^{\sigma} \|^{-1} u_n^{\sigma}~~.\]
Since $\Delta_G^{-1/2} \le c^{n+1}$ on the support of $f_n$,
we have $\| \Delta_G^{-1/2} f_n \| \le c^{n+1} \| f_n \|$. 
On the other hand, $ \| \Delta_G^{-1/2} v_n^{\sigma} \|^{-1}
\le c^{-n}$, hence $\| A(\sigma)^* C_{\sigma} f\| \le c \| f \|$, which means that 
$[ A(\sigma)^* C_{\sigma} ]$ exists. Thus $A(\sigma)$ is as desired.

Let us next address the measurability requirement: With respect to
direct integrals, we use the terminology of \cite{Fo}. Clearly it
is sufficient to show that the orthonormal basis $(u_n^{\sigma})_{n \in \NN}$
and the images $(\| \Delta_G^{-1/2} v_n^{\sigma} \|^{-1} v_n^{\sigma})_{n \in
\NN}$ can be
chosen measurably, that is, for each $n \in \NN$, $(u_n^{\sigma})_{\sigma
\in V}$ and $(\| \Delta_G^{-1/2} v_n^{\sigma} \|^{-1} v_n^{\sigma})_{\sigma
\in V}$ are measurable vector fields.

First of all, since $\tau(V)$ is standard, there exists a measurable
realization of the family $(\tau(\sigma), {\mathcal H}_{\tau(\sigma)})_{\sigma
\in V}$ \cite[Theorem 10.2]{Ma}. This means that there exist vector
fields 
$((e_k^{\sigma})_{\sigma \in V})_{k \in
\NN}$, such that for almost all $\sigma \in V$,
$(e_k^{\sigma})_{k \in \NN}$ is total in 
${\mathcal H}_{\tau(\sigma)}$, and the mappings $\sigma \mapsto 
\langle e_k^{\sigma},e_m^{\sigma} \rangle$ as well as
$\sigma \mapsto \langle e_k^{\sigma}, \tau(\sigma) (g_0) e_m^{\sigma} \rangle$
are measurable, for all $n,m \in \NN$ and all $g_0 \in G_0$.
The measurability of any
vector field is equivalent to the measurability of its scalar products
with the $e_k^{\sigma}$. Without loss of generality, we may assume in
addition that, for $\nu_G$-almost every
$\sigma$, the first ${\rm dim}({\mathcal H}_{\tau(\sigma)})$ vectors
are an orthonormal basis of ${\mathcal H}_{\tau(\sigma)}$ \cite[7.29, 7.30]{Fo}. 
Let $(a_n)_{n \in \NN}$ be any orthonormal basis
of ${\rm L}^2(\Gamma)$, then $(a_n e_k^{\sigma})_{n,k}$ is an orthonormal
basis of $ {\mathcal H}_{\sigma} = {\rm L}^2(\Gamma, d\mu_\Gamma; 
{\mathcal H}_{\tau(\sigma)})$, except for the zero vectors belonging to
the indices $k > {\rm dim}({\mathcal H}_{\tau(\sigma)})$. 
This family makes $({\mathcal H}_{\sigma})_{\sigma \in V}$ a measurable
field of Hilbert spaces, and we are going to ensure
measurability of the operator fields with respect to this structure.

In order to construct the measurable family $(u_n^{\sigma})_{n,\sigma}$ of 
orthonormal bases, we proceed as follows:
First note that the sets $V_\ell := \{\sigma: {\rm dim}
({\mathcal H}_{\tau(\sigma)}) = \ell \}$, for $\ell \in \NN \cup \{ \infty \}$, 
are Borel sets \cite[Theorem 8.7]{Ma}. On each $V_\ell$, pick a fixed bijection
$s_\ell: \NN \times \{ 1, \ldots, \ell \} \to \NN$ (where $\{ 1, \ldots,
\infty \} := \NN$). Then  letting
$u_{s_\ell(n,k)}^{\sigma} := a_n e_k^{\sigma}$, for $\sigma \in
V_\ell$, removes the zero
vectors. Moreover, on each $V_\ell$, the measurability is easily checked,
and this is sufficient.

For the construction of the $v_n$ we pick any family $(b_n)_{n \in \NN} 
\subset {\rm L}^2(\Gamma)$ of unit vectors, such that $b_n$ is supported in
$S_n$. Moreover, let $(\xi^{\sigma})_{\sigma \in V}$ be a measurable field
of unit vectors $\xi^{\sigma} \in {\mathcal H}_{\tau(\sigma)}$, and define
$v_n^{\sigma} = b_n \xi^{\sigma}$.
Then
\[ \sigma \mapsto \langle \| \Delta_G^{-1/2} a_n \|^{-1}
v_n^{\sigma} , a_n e_k^{\sigma} \rangle =  \| \Delta_G^{-1/2} a_n \|^{-1}
 \langle b_n, a_n \rangle \langle \xi^{\sigma}, e_k^{\sigma}\rangle ~\]
is measurable by the choice of the $\xi^{\sigma}$. Thus we can
construct the operator field in a measurable way.

Finally, let us provide for square-integrability. For this
purpose we observe that we may assume the constant $c$ picked above 
to be $\ge 2$, and then $\| A(\sigma) \|_2^2 < 2$. Moreover, if we shift
the construction in the sense that $u_n^{\sigma}\mapsto\|
\Delta_G^{-1/2} v_{n+k}^{\sigma} \|^{-1} v_{n+k}^{\sigma}$, for
$k > 0$, we obtain $\| A(\sigma) \|_2^2 < 2^{1-k}$, while preserving all
the other properties of $A(\sigma)$. With this in mind, we can easily
modify the construction to obtain an element of ${\mathcal B}_2^{\oplus}$:
Since $G$ is separable, $\nu_G$ is $\sigma$-finite, i.e., 
$\hat{G} = \bigcup_{n \in \NN} \Sigma_n$ with the
$\Sigma_n$ pairwise disjoint and $\nu_G(\Sigma_n) < \infty$.
Shifting on $\Sigma_n$ by $k_n \in N$ with $\nu_G(\Sigma_n) 2^{-k_n}
< 2^{-n}$ ensures square-integrability without destroying measurability.
(The latter is obvious on $\Sigma_n$.)
Hence we are done, noting that the shifting argument also yields the
following fact which sharpens the contrast to the unimodular case,
where for a given representation the length of admissible vectors is fixed.

\begin{cor}
 There exist admissible vectors with arbitrarily small or big norm.
\end{cor}

\section{Concluding remarks}

The main purpose of this paper was to establish a link between
generalized wavelet transforms and Plancherel theory, and to 
demonstrate the power of the approach via the existence theorem
for admissible vectors. It is quite obvious that for a concrete
situation the abstract approach does not automatically give
access to the construction of admissible vectors. 
It can be difficult to establish the fact that a given
representation is contained in the regular representation, and
even when that is known, obtaining explicit knowledge of the decomposition of 
the representation under the Plancherel transform might turn out to
be another serious obstacle. Nevertheless, the approach yields a new 
perspective on the construction of wavelet transforms, and there are concrete 
cases where the link to Plancherel theory can be established. 
A family of examples are the semidirect products $G = \RR^n \rtimes H$,
and their quasi-regular representations already mentioned
in the introduction, and studied, with increasing generality,
in the papers \cite{Mu,IsKl,KlSt,BeTa,Fu,LWWW}. In these cases both
the decomposition of the quasi-regular
representation and the Plancherel theory of $G$ can be computed
explicitly, and
it is instructive to view the various admissibility conditions
derived for those groups in the light of the abstract approach.
We intend to discuss this point in more detail in a further publication.

\section{Acknowledgements}

The author whishes to thank the Laboratoire d'Analyse, Topologie et
Probabilit\'{e}s
in Marseille, and the Department of Mathematics and Statistics of
Concordia University, Montreal, for their hospitality. He has
profited from stimulating discussions with S.T. Ali.
The support
by the Deutsche Forschungsgemeinschaft (DFG) under the contract
Fu 402/1 is also acknowledged.


\begin{thebibliography}{99}
 \bibitem{AlAnGa} {S. T. Ali, J-P. Antoine and J-P. Gazeau,
  {\em Coherent States, Wavelets and Their Generalizations},
  Springer-Verlag,  New York, 2000.}
  \bibitem{AlFuKr} {S.T. Ali, H. F\"uhr and A. Krasowska, {\em
  Wavelet transforms, Wigner functions and Plancherel inversion,}
  in preparation.}
  \bibitem{ArLu}{D. Arnal and J. Ludwig, {\em Q.U.P. and Paley-Wiener
  property of unimodular, especially nilpotent, Lie groups,}
  Proc. Amer. Math. Soc. {\bf 125}, (1997), 1071-1080.}
 \bibitem{BeTa}{D. Bernier and K. Taylor, {\em Wavelets from 
  square-integrable representations,} SIAM J. Math. Anal. {\bf 27} (1996), 
  594-608.}
 \bibitem{Bo} {G. Bohnke, {\em Treillis d'ondelettes aux groupes de
 Lorentz,} Annales de l'Institut Henri Poincar\'{e} {\bf 54} (1991),
 245-259.}
  \bibitem{Ca} {A.L. Carey, ``Group representations in reproducing 
kernel Hilberts spaces,'' Reports in Math. Phys. {\bf 14} (1978), 247--259}
 \bibitem{Di} {J. Dixmier, {\em $C^{\ast}$-Algebras,}
  North Holland, Amsterdam, 1977.}
  \bibitem{DuMo} {M. Duflo and C.C. Moore, {\em On the regular 
  representation of a nonunimodular locally compact group,} J. Funct.
  Anal. {\bf 21} (1976), 209-243.}
 \bibitem{Fo} {G.B. Folland, {\em A Course in Abstract Harmonic Analysis,}
  CRC Press, Boca Raton, 1995.}  
 \bibitem{Fu} {H. F\"uhr, {\em Wavelet frames and admissibility in higher
 dimensions,} J. Math. Phys. {\bf 37} (1996), 6353-6366.}
 \bibitem{GrMoPa}{A. Grossmann, J. Morlet and T. Paul, 
  {\em Transforms associated to square integrable group representations I:
  General Results,} J. Math. Phys. {\bf 26} (1985), 2473-2479.} 
 \bibitem{IsKl}{C.J. Isham and J.R. Klauder, {\em Coherent states for
  $n$-dimensional Euclidean groups $E(n)$ and their application,}
  J. Math. Phys. {\bf 32} (1991), 607-620.}
 \bibitem{Ji} {Q. Jiang, {\em Wavelet transform and orthogonal decomposition
 of ${\rm L}^2$ space on the Cartan domain $BDI(q=2)$,}
 Trans. Am. Math. Soc. {\bf 349} (1997), 2049-2068.}
 \bibitem{KlSt} {J.R. Klauder and R.F. Streater, {\em A wavelet transform
  for the Poincar\'{e} group,} J. Math. Phys. {\bf 32} (1991), 1609-1611.}
 \bibitem{LWWW} {R.S. Laugesen, N. Weaver, G. Weiss and E.N. Wilson, 
 {\em A generalized Calder\'on reproducing formula and its associated
 continuous wavelets,} in preparation.}
 \bibitem{Li} {R.L. Lipsman: {\em Non-abelian Fourier analysis.}
 Bull. Sci. Math. {\bf 98} (1974), 209-233.}
 \bibitem{Ma} {G. Mackey, {\em Borel structure in groups and their duals,}
 Trans. Amer. Math. Soc. {\bf 85} (1957), 134-165.}
 \bibitem{Mu} {R. Murenzi, {\em Ondelettes multidimensionelles et 
  application \`{a} l'analyse d'images,} Th\`{e}se, 
  Universit\'{e} Catholique de Louvain, Louvain-La-Neuve, 1990.}
\bibitem{Ta} { N.~Tatsuuma, {\em Plancherel formula for non-unimodular
 locally compact groups,} J. Math. Kyoto Univ. {\bf 12} (1972), 179-261.}
\end{thebibliography}
\end{document}